# Real-Time Optimal Power Flow under Wind Energy Penetration-Part I: Approach

Erfan Mohagheghi, *Student Member*, *IEEE*, Aouss Gabash, *Member*, *IEEE*, Pu Li
Department of Simulation and Optimal Processes
Institute of Automation and Systems Engineering
Ilmenau University of Technology
Ilmenau, Germany
erfan.mohagheghi@tu-ilmenau.de, aouss.gabash@tu-ilmenau.de, pu.li@tu-ilmenau.de

*Abstract*—Real-time optimal power flow (RT-OPF) under wind energy penetration is highly desired but extremely difficult to realize. This is basically due to the conflict between the fast changes in wind power generation and the slow response from the optimization computation. This paper (Part I) presents a prediction-updating approach to address this challenge. We consider essential scenarios around forecasted data of wind power that would probably happen during the computation time required for solving a large-scale complex optimal power flow problem. Parallel computing is used to solve the individual OPF problems corresponding to these scenarios. This provides for the forecasted time horizon probable reference operations in the form of a lookup-table. One of these operations will be selected based on the actual wind power and realized to the grid for the current time interval, thus leading to a RT-OPF framework. The proposed approach is implemented in Part II of this paper using a 41-bus medium-voltage distribution network as a case study.

## NOMENCLATURE

**Set**

| | |
|---|---|
| $sl$ | Set of wind power levels, i.e., {1, 2, 3}. |

**Functions**

| | |
|---|---|
| $f$ | Objective function. |
| $g$ | Equality equations. |
| $h$ | Inequality equations. |

**Parameters**

| | |
|---|---|
| $P_d(i,t_{120s})$ | Active power demand at bus $i$ in prediction horizon $t_{120s}$. |
| $P_w$ | Active power of a wind station (WS). |
| $P_{W.r}(i)$ | Rated installed wind power at bus $i$. |
| $P_{w.A}(i,t_{20s})$ | Actual wind power of WS at bus $i$ in update interval $t_{20s}$. |
| $P_{w.H.\sigma}(i,t_{120s})$ | Wind power higher than forecasted of WS at bus $i$ in prediction horizon $t_{120s}$ for level $\sigma \in sl$. |
| $P_{w.L.\sigma}(i,t_{120s})$ | Wind power lower than forecasted of WS at bus $i$ in prediction horizon $t_{120s}$ for level $\sigma \in sl$. |
| $P_{w.M}(i,t_{120s})$ | Mean (forecasted) wind power of WS at bus $i$ in prediction horizon $t_{120s}$. |
| $Q_d(i,t_{120s})$ | Reactive power demand at bus $i$ in prediction horizon $t_{120s}$. |
| $t_{120s}$ | Prediction horizon, i.e., 120 seconds. |
| $t_{112s}$ | Reserved time for computing OPF problems, i.e., 112 seconds. |
| $t_{20s}$ | Update interval, i.e., 20 seconds. |
| $t_{2s}$ | Reserved time for data management, i.e., 2 seconds. |
| $u^{max}$ | Upper limit on control variables. |
| $u^{min}$ | Lower limit on control variables. |
| $x^{max}$ | Upper limit on state variables. |
| $x^{min}$ | Upper limit on state variables. |
| $\Delta P_\sigma(i)$ | Wind power deviation at bus $i$ for level $\sigma \in sl$. |

**Variables**

| | |
|---|---|
| $P_S$ | Active power imported from slack bus. |
| $Q_S$ | Reactive power imported from slack bus. |
| $u$ | Vector of control variables. |
| $x$ | Vector of state variables. |
| $\beta_{c.w}(i,t_{120s})$ | Curtailment factor of wind power for WS at bus $i$ in prediction horizon $t_{120s}$. |
| $\beta_{c.w}(i,t_{20s})$ | Curtailment factor of wind power for WS at bus $i$ in update interval $t_{20s}$. |

## I. INTRODUCTION

The dramatic increase of renewable energy penetration leads to a significant challenge in the operation of energy distribution networks (DNs). In particular, the wind power generation is intermittent, i.e., the DN operator has to fast update the operation strategy correspondingly. It is highly desired to carry out this task by an online optimization. However, the optimization problem for this task is usually so large and complicated that it takes the computation time which is much higher than required for compensating the fast change of the wind power generation. Even by using advanced optimization algorithms combined with modern computation facilities, the computation time is not short enough to achieve this target. Therefore, efficient and reliable computation approaches need to be further developed for the

This work is supported by the Carl-Zeiss-Stiftung.

The final version of this paper has been published in the proceeding of 2016 IEEE 16th International Conference on Environment and Electrical Engineering (EEEIC). ©2016 IEEE. DOI: 10.1109/EEEIC.2016.7555464

implementation of real-time optimal power flow (RT-OPF) under wind energy penetration.

Introduced by Carpentier in 1962 [1], optimal power flow (OPF) has been widely used for operation planning of power networks [2]-[4], where no renewable energy generation (REG) was considered in the network. In contrast, OPF with REG was taken into account in [5] and recently, an active-reactive OPF in active DNs was introduced in [6], [7]. One of the remarkable abilities of the methods in [6], [7] is its capability to ensure feasible and optimal solutions even with a high penetration of wind and solar energy. This is achieved by using a curtailment strategy [8]-[10], by which a part of wind power generation is curtailed in order to satisfy system constraints. However, OPF usually leads to a large-scale complex optimization problem as the objective function becomes more and more complicated, especially if the contribution of costs and/or revenues is modelled in more detail [6]. As a result, the solution of such an optimization problem remains a challenging task, in particular for the realization of RT-OPF.

In most above mentioned studies on OPF, the forecasted values of renewable energy penetration and demand were used for a prediction horizon, whereas other studies have been made to improve the quality or accuracy of forecasting [11], [12]. However, there always exist discrepancies between the forecasted and the actual values, especially for wind power generation [13]. Such discrepancies will lead to constraint violations, i.e., the OPF results obtained in this way cannot be used for the practical operation. Therefore, fluctuations or deviations from the forecasted values of the wind energy penetration in the prediction horizon, in which the operation is to be planned, should be considered in the OPF computation.

In general, there are two approaches to address this problem. First, a stochastic optimization problem can be formulated in which the demand and/or REG are treated as uncertain parameters with a known stochastic distribution. This approach is known as probabilistic OPF [14]. A promising method here is to formulate and solve a chance constrained optimization problem [15], [16] which can provide optimal and reliable operation strategies under uncertainty. However, these probabilistic methods are in effect offline approaches, i.e., the computation for solving the optimization problem is highly time-consuming [15], especially when considering real-time implementation.

The second approach to handling the discrepancies between the forecasted and actual values of demand and/or REG is to perform a prediction-updating strategy. The controls will be at first derived based on the forecasted values for a future time horizon and then updated based on the actual values when available. The concept of RT-OPF was early demonstrated in [17] by solving a quadratic optimization problem. In [18], a real-time optimal reactive power flow was executed for several selected loading conditions, where the frequency of executions could vary between 15 minutes up to four hours. It is to note that [17], [18] did not consider any REG in the network. The authors in [19] proposed an energy management system using neural networks trained on some scenarios of the uncertain demand and REG. In [19], a RT-OPF scheme for a 23-bus radial DN with two wind turbine generators was used to verify the method, where the time for updating the solutions was less than three minutes.

Recently, the effect of curtailment levels of two power flow management approaches, i.e., the constraint satisfaction [20] and OPF method [21], was evaluated in [10]. Both methods were based on heuristic rules to determine the curtailment and thus a closed-loop implementation could be realized. Since the rule-based methods react fast, a time horizon was not required to predict the future wind power, while in contrast a risk-based AC-OPF approach was proposed as an online framework in [22] to handle the uncertain wind power in a future time horizon. The resulting optimal set-points for wind power curtailment were only to ensure the safety of the operation, i.e., no economic aspects were considered in [22].

Based on the above review of the recent studies, real-time optimization is desired in which not only technical but also economic (e.g., energy prices, cost of grid energy losses, feed-in tariffs and reverse power flow [23], [24]) issues need to be considered for system operators to achieve both reliable and optimal operations. In particular, due to the computation time needed for solving the OPF problem, the future wind power generation and/or load in a prediction time horizon has to be predicted. Therefore, it is aimed in this paper to develop a new techno-economic RT-OPF framework consisting of following innovative aspects:

- Minimizing the costs of wind power curtailment considering technical and economic issues simultaneously;
- Wind power curtailment will be optimized in real-time under intermittent wind power (IWP) penetration;
- Considering both active and reactive energy prices while minimizing the costs of grid active energy losses;
- Realizing optimal operation strategies to the grid in a practically desired sampling time (every 20 seconds);
- Reducing the issue of computational time by utilizing parallel computing.

The basic idea of our approach is to use forecasted wind power as an expected profile and consider its probable scenarios around these profiles in a prediction horizon to solve the optimal operation problem. The resulting operation strategy in the first time interval will be realized to the power system. The planned strategy will be updated based on the currently measured wind power. This leads to a prediction-updating framework in a way of a moving horizon scheme.

The remainder of the paper is organized as follows. Section II describes the problem of OPF considering IWP and energy prices. The new techno-economic RT-OPF framework to deal with IWP is proposed in section III. The paper is summarized and discussed in section IV.

II. PROBLEM DESCRIPTION

*A. OPF under IWP*

In [6], an OPF problem was formulated with forecasted input wind power and demand profiles in the prediction horizon for operation planning. Such a deterministic OPF problem can reduce the computational effort on the one hand, but on the other hand it may fail to handle the impreciseness of the forecasted inputs. Therefore, our aim is to overcome this problem by focusing on detailed consideration of IWP. We consider a medium-voltage (MV) network with two wind stations (WSs) as shown in Fig. 1. Here, active $P_S$ and reactive $Q_S$ power at slack bus are allowed to be either positive (imported from a high-voltage (HV) network) or zero (no reverse power flow to the HV network) to avoid any possible generation rejections [23], i.e., active and reactive power can be imported (in the case of low wind power), but not exported (in the case of high wind power) [24]. In addition, a low-voltage (LV) network consumes active and reactive power as shown in Fig. 1. Due to system constraints, a wind power curtailment factor ($0 \leq \beta_{c.w} \leq 1$) at each WS is used as a control variable [6], where $\beta_{c.w} = 1$ when no curtailment and $\beta_{c.w} < 1$ otherwise.

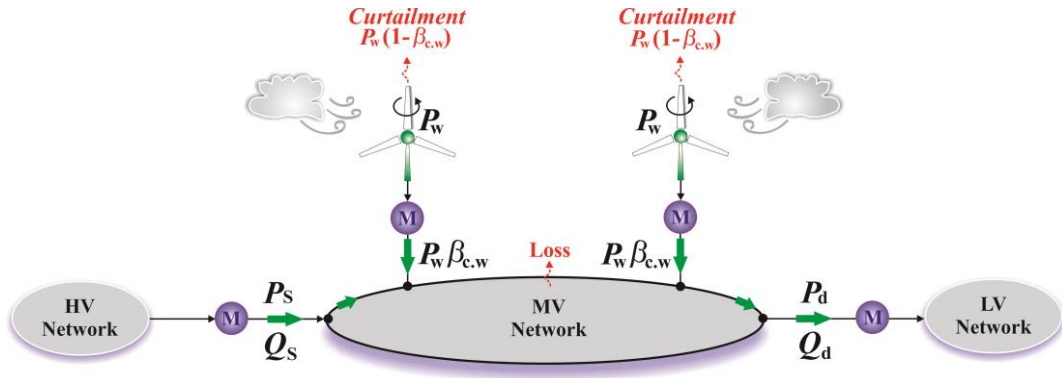

Fig. 1. Illustration of a meter-based method for charging and remunerating different entities connected to a power system [24].

In addition, we assume that the two WSs have the same installed capacities. The wind power $P_w$ generated from each WS can be forecasted from the expected wind speed [6] at every forecasting time horizon, i.e., a prediction horizon ($t_{120s}$), see Fig. 2. It means that the forecasted wind power is considered as a constant in this future time horizon and will be updated at each $t_{120s}$. Note that $P_w$ from the two WSs can be different even with a short distance between them. From another perspective, the actual wind power (AWP) in the time horizon $t_{120s}$ can be different from the forecasted value of $P_w$. Thus we consider several essential scenarios around the forecasted value to describe the variations in wind power generation for each WS, as shown in Fig. 2. These scenarios represent highly probable events of the IWP and can be chosen based on the forecasted data. Here, $P_{w.M}$ stands for the mean (the forecasted wind power), $P_{w.H.\sigma}$ for the higher-side (values higher than forecasted) and $P_{w.L.\sigma}$ for the lower-side (values lower than forecasted), see Fig. 2. Here, three levels on each side are considered. It is also to note that the wind power deviation, (the width of each level) should be ($0 \leq \Delta P_\sigma(i) \leq P_{w.r}(i)$), where $\Delta P_\sigma(i)$ is the wind power deviation at bus $i$ for level $\sigma$ and $P_{w.r}(i)$ is the rated power of WS at bus $i$. Four time horizons are illustrated in Fig. 2, where the black line is the expected (forecasted) profile while the red and green dashed lines are the boundaries of the levels on the higher-side as well as on the lower-side, respectively. The black dots represent the actually generated values of wind power $P_{w.A}$ for a WS.

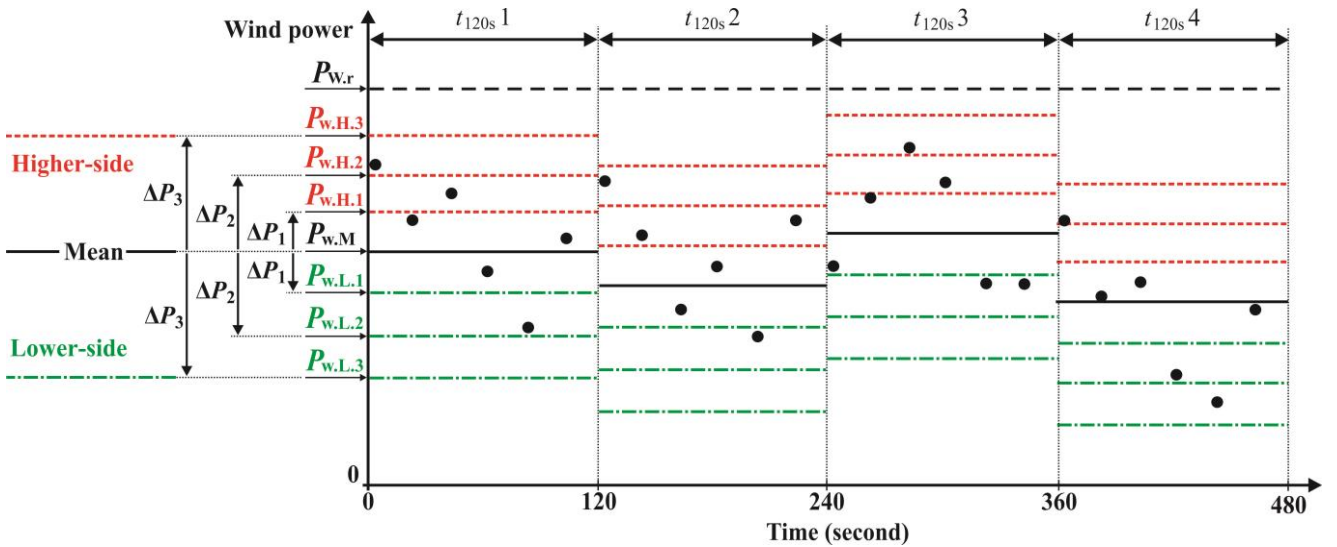

Fig. 2. General illustration of wind power for a WS in four prediction horizons.

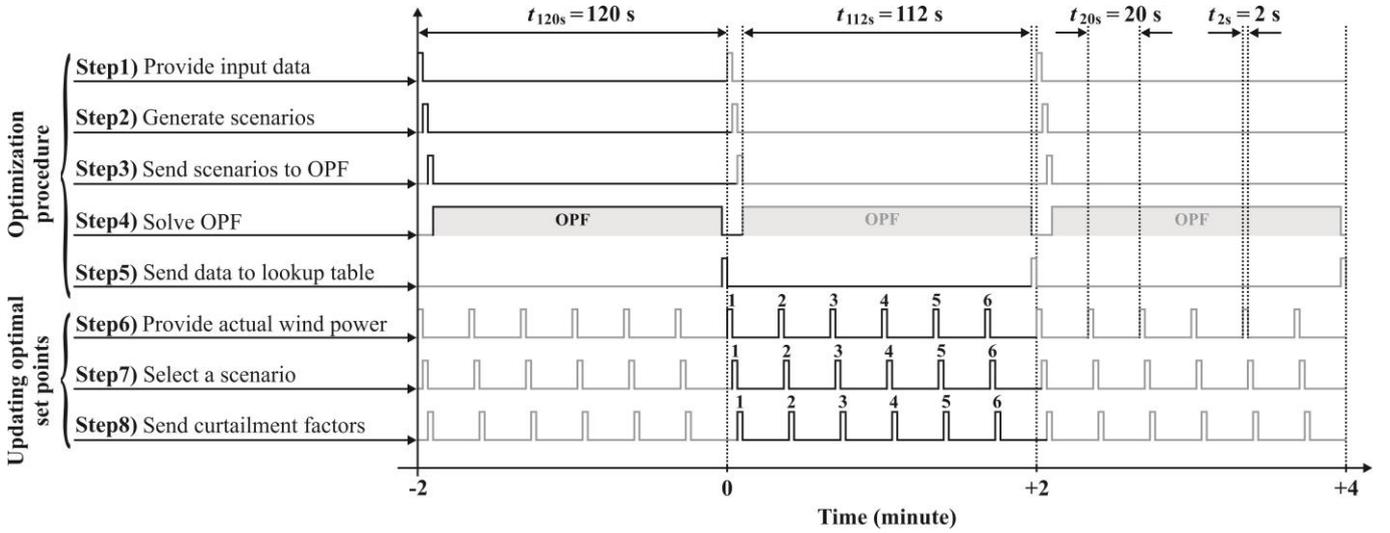

Fig. 3. Time allocation for the computational tasks of the 8 steps. Here, 3 prediction horizons are shown.

*B. Active and Reactive Energy at Slack Bus*

Our objective of OPF for the grid shown in Fig. 1 is to maximize the total revenue from the wind power (i.e., to minimize wind power curtailment) and meanwhile to minimize the total costs of the active energy losses in the MV network as well as the total costs of the active and reactive energy imported from the HV network. The active and reactive energy at the slack bus in this work are defined as follows:

- *Active Energy*: The forward active energy from the HV network to the MV network is to be minimized based on an active energy price model, while the reverse active energy is not allowed in order to avoid any possible active power rejection from the HV network [23].

- *Reactive Energy*: The forward reactive energy from the HV network to the MV network is to be minimized based on a reactive energy price model, while the reverse reactive energy is not allowed in order to avoid any possible charge for reactive energy [25].

As a result, the system does not export any active and reactive power to the upstream network. It is worth mentioning that, so far, there has been no definitive solution on how reverse power flow can be realized, as recently discussed in [26], [27].

*C. Problem Formulation*

In this work, the OPF problem under IWP is formulated as a nonlinear optimization problem. The general formulation of the problem is as follows:

$$\max_{u} f(x,u) \quad (1)$$

where $f$ is the objective function to be maximized which is a trade-off between the total revenue from active wind power injected to the MV network, the total costs of active energy losses in the gird, and the costs of the active and reactive energy at the slack bus, detailed mathematical formulations are given in Part II of our paper [28]. Here, $u$ represents the vector of control variables (curtailment factors of WSs) and $x$ is the vector of state variables (real and imaginary parts of bus complex voltage at PQ buses and active and reactive power at the slack bus). The objective function $f$ is subject to equality $g$ and inequality $h$ equations as follows

$$g(x,u) = 0 \quad (2)$$

$$h(x,u) \leq h^{\max} \quad (3)$$

$$x^{\min} \leq x \leq x^{\max} \quad (4)$$

$$u^{\min} \leq u \leq u^{\max}. \quad (5)$$

III. PROPOSED RT-OPF FRAMEWORK

The proposed RT-OPF framework in this paper consists of five steps for the prediction and three steps for the updating, as shown in Fig. 3. The input data for this framework include forecasted wind power generation of the WSs at different locations in the MV network for the prediction horizon $t_{120s}$. Here, $t_{120s}$ is taken 120 seconds which is divided into 6 sub-intervals each of which (denoted as $t_{20s}$) is 20 seconds, as shown in Figs. 2 and 3. The steps shown in Fig. 3 are explained below:

**Step 1)** Provide active and reactive power demand, $P_d(i,t_{120s})$ and $Q_d(i,t_{120s})$, as inputs to the OPF. At the same time, forecasted wind power $P_{w.M}(i,t_{120s})$ for the current prediction horizon is provided.

**Step 2)** Generate scenarios for the higher-side and lower-side levels, i.e., $P_{w.H.\sigma}(i,t_{120s})$ and $P_{w.L.\sigma}(i,t_{120s})$, for the prediction horizon using

$$P_{w.H.\sigma}\left(i,t_{120s}\right) = P_{w.M}\left(i,t_{120s}\right) + \Delta P_{\sigma}(i), \qquad \sigma \in sl \qquad (6)$$

$$P_{w.L.\sigma}\left(i,t_{120s}\right) = P_{w.M}\left(i,t_{120s}\right) - \Delta P_{\sigma}(i), \qquad \sigma \in sl \qquad (7)$$

where $\Delta P_3(i) = 1.5\Delta P_2(i) = 3\Delta P_1(i)$, see Fig. 2. The total number of scenarios is the number of the combination of the seven levels of each WS. Note that the scenarios generated should be between 0 and the rated power of the WS, i.e., if $P_{w.H.\sigma}(i,t_{120s}) > P_{W.r}(i)$, then $P_{w.H.\sigma}(i,t_{120s}) = P_{W.r}(i)$ and if $P_{w.L.\sigma}(i,t_{120s}) < 0$, then $P_{w.L.\sigma}(i,t_{120s}) = 0$.

**Step 3)** Deliver the calculated values of the wind power scenarios (i.e., $P_{w.L.\sigma}(i,t_{120s})$, $P_{w.M}(i,t_{120s})$, and $P_{w.H.\sigma}(i,t_{120s})$) as inputs to the optimization computation.

**Step 4)** Solve the OPF for all scenarios (49 scenarios for two WSs and seven levels for each WS). Since each scenario is independent, parallel computing is used here to obtain the results within the reserved time $t_{112s}$ for the computation.

**Step 5)** Form the scenarios and corresponding OPF results including $\beta_{c.w}(i,t_{120s})$ into a lookup table. The lookup table will be updated for every prediction horizon $t_{120s}$.

**Step 6)** The actual wind power $P_{w.A}(i,t_{20s})$ (supposed to be available each time interval $t_{20s}$) is provided.

**Step 7)** Select one of the scenarios and the corresponding $\beta_{c.w}(i,t_{120s})$ based on the following rule. If $P_{w.A}(i,t_{20s})$ is not equal to $P_{w.L.\sigma}(i,t_{120s})$, $P_{w.M}(i,t_{120s})$, or $P_{w.H.\sigma}(i,t_{120s})$, then consider it to be as *the nearest higher level*, since higher wind power corresponds with *higher risk* and lower controls correspond with *lower risk*.

**Step 8)** Realize the values of the optimal controls $\beta_{c.w}(i,t_{20s})$ corresponding to the selected scenario to the WSs. It means that the computed optimal amount of wind power is penetrated to the MV network in the current time interval $t_{20s}$. The Steps 6, 7, and 8 will be repeated from one interval $t_{20s}$ to the next (6 times) till to the end of the prediction horizon $t_{120s}$, with which one cycle of the RT-OPF is completed and the next cycle starts from Step 1. In Fig. 3, the cycle is shown by black lines distinguished from the gray lines.

It should be noted that the decision made in Step 7 has a certain degree of conservatism to ensure the feasibility of the determined curtailment factors to be realized to the grid. This means that, guaranteeing the feasibility as a higher priority, we have to sacrifice some amount of generated wind power for safe operation (also see [29]). This sacrificed amount can be reduced if more levels on the higher- and lower-side are defined, but then the number of scenarios will be increased correspondingly.

To ensure a unified time horizon $t_{120s}$ in the computation sequences in Fig. 3, we need to consider the time allocations for the different tasks. In Fig. 3, $t_{2s}$ denotes the reserved time for data management (in this paper 2 seconds). For Steps 1, 3, 5, 6, and 8, $t_{2s}$ means the communication time for sending/receiving data. In Steps 2 and 7, $t_{2s}$ means the time for processing data after receiving the corresponding inputs and $t_{112s}$ is the time reserved for the processors to solve the OPF problems. Here, $t_{112s} = 112$ seconds to ensure that none of the processors exceed this limit, since the OPF computation takes the largest part of the time horizon. Note that the reserved time for computing the OPF problems can be different due to the network size and complexity.

The proposed RT-OPF framework is implemented on a medium-voltage DN, which is reported in Part II [28] of this paper.

## IV. SUMMARY

In this paper, a real-time optimal power flow (RT-OPF) framework is proposed for active distribution networks with intermittent wind power penetration. It ensures that the OPF operation strategy will be updated in a short sampling time (i.e., every 20 seconds) based on the real wind power values. This is achieved by solving the optimization problems for a prediction horizon (i.e., 120 seconds) corresponding to probable scenarios around the forecasted wind power. The proposed approach provides, in real-time, not only optimal but also feasible solutions for the grid operation. The scenario-based individual optimization problems are solved in parallel to guarantee the computation in the reserved time.